\documentclass{amsart}[10]
\usepackage{amsmath}
\usepackage{amsfonts}
\usepackage{amssymb}
\usepackage{verbatim}

\usepackage[title]{}

\newcommand{\PP}{{\bf P}}

\newcommand{\Q}{{\mathbb Q}}
\newcommand{\Z}{{\mathbb Z}}
\newcommand{\N}{{\mathbb N}}
\newcommand{\C}{{\mathbb C}}
\newcommand{\A}{\mathcal{C}}
\newcommand{\Ac}{\mathbb{A}}

\newcommand{\bl}{\subsection{Lemma}}
\newcommand{\el}{ \medskip  }
\newcommand{\bt}{\subsection{Theorem}}
\newcommand{\et}{\medskip  }
\newcommand{\bp}{\subsection{Proposition}}
\newcommand{\ep}{\medskip  }
\newcommand{\bc}{\subsection{Corollary}}
\newcommand{\ec}{\medskip  }
\newcommand{\df}{\subsection{Definition}}
\newcommand{\edf}{\medskip  }
\newcommand{\lb}{\label}

\newcommand{\be}{\begin{equation}}
\newcommand{\ee}{\end{equation} }

\newcommand{\pr}{{\rm pr\, }}
\newcommand{\prs}{{\rm pr}}
\newcommand{\cl}{{\rm cl}}

\newcommand{\ark}{{\rm ark}}

\newcommand{\dd}{\partial}

\newcommand{\ra}{\rangle}
\newcommand{\la}{\langle}

\newcommand{\inv}{^{-1}}

\newcommand{\M}{{\rm M}}

\newcommand{\subs}{\!\subseteq\!}
\newcommand{\sbo}{\subseteq_{op}}
\newcommand{\sba}{\subseteq_{an}}
\newcommand{\sups}{\!\supseteq\!}
\newcommand{\smin}{\setminus}
\newcommand{\stmn}{\setminus}

\newcommand{\pf}{{\bf Proof. }}

\newcommand{\concat}{{{}^{\frown}}}

\newcommand{\rrem}{\subsection{Remark}}
\newcommand{\erem}{\medskip}

\newcommand{\F}{K}

\newcommand{\exf}{{\rm ex}}

\newcommand{\sbc}{\subseteq_{\rm cl}}

\newcommand{\dom}{{\rm Dom\,}}

\newcommand{\bex}{\begin{ex}}
\newcommand{\eex}{\end{ex}}

\title{Analytic Zariski structures and 
non-elementary categoricity}

\author{Boris Zilber}
\address{University of Oxford}

\date{12 January, 2016}
\begin{document}

\maketitle
\begin{abstract}{We study  analytic Zariski structures from the point of view of
non-elementary model theory. We show how to associate an abstract elementary 
class with a one-dimensional analytic Zariski structure and prove that the 
class is stable, quasi-minimal and homogeneous over models.
We also demonstrate how Hrushovski's predimension arises in this 
general context as a natural geometric notion and use it as one of our main tools. 
  }
\end{abstract}
The notion of an  analytic Zariski structure was introduced in \cite{[1]} by
the author and N.Peatfield in a form slightly different from the one
presented here and then in \cite{[2]}, Ch.6  in the current form.  Analytic Zariski generalises the previously known
notion of a Zariski structure (see \cite{[2]} for one-dimensional case and \cite{[3]}, \cite{[4]} for the general definition) mainly by 
dropping the
requirement of Noetherianity and weakening the assumptions on the
projections of closed sets.  We remark that in the broad setting, it is appropriate to consider the notion of a Zariski structure as belonging to positive model theory in the sense of \cite{Be}.

In \cite{[1]} we assumed that the Zariski structure is compact (or compactifiable), 
here we drop this assumption, which may be too restrictive in applications. 

The class of analytic Zariski structures is much broader and geometrically richer than the class of Noetherian Zariski structures. 
The main examples come from two sources:

(i) structures which are constructed in terms of complex analytic
functions and relations;

(ii) ``new stable structures'' introduced by Hrushovski's construction;
 in many cases these objects exhibit properties similar
to those of class (i).

However, although there are concrete examples for both (i) and (ii), 
in many  cases we lack the technology to prove  that the structure 
is analytic Zariski. In particular, despite some attempts the 
conjecture that $\C_{\exp}$ is analytic Zariski, 
assuming it satisfies axioms of 
pseudo-exponentiation (see \cite{[14]}),
  is still open. \\


The aim of this paper is to carry out a model-theoretic analysis of  $\M$ in the appropriate language. 
Recall that if $\M$ is a Noetherian Zariski structure the relevant key model-theoretic result states that its first-order theory allows elimination of quantifiers and is $\omega$-stable of finite Morley rank. In particular, it is strongly minimal (and so uncountably categorical) if $\dim M=1$ and $\M$ is irreducible.

For analytic Zariski 1-dimensional $\M$ we carry out a model theoretic study it in the spirit of the theory of {\em abstract
elementary classes}. We start by introducing a suitable countable fragment
of the family of basic Zariski relations and a correspondent substructure of constants over which all the further analysis is carried out.
Then we proceed to the analysis of the notion of dimension of Zariski closed sets and 
define more delicate  notions of the {\em predimension} and {\em dimension}
of a tuple in $\M$. 
In fact by doing this  we reinterpret dimensions which are present in 
every analytic structure in terms familiar to many from Hrushovski's construction, thus establishing once again conceptual links between classes (i) and (ii).

Our main results are proved under assumption that $\M$ is one-dimensional
(as an analytic  Zariski structure) and irreducible. No assumption on
presmoothness is needed. We prove for such an $\M,$ in the terminology of 
 \cite{[13]}:

(1) {\em $\M$ is a quasi-minimal pregeometry structure with regards to a closure operator $\cl$ associated
with the predimension;}

(2) {\em $\M$ has quantifier-elimination to the level of $\exists$-formulas in the following sense:
every two tuples which are (first-order) $\exists$-equivalent over a countable submodel, are
$L_{\infty,\omega}$ equivalent;}

(3) {\em The abstract elementary class associated with $\M$ is categorical in uncountable cardinals and is excellent.}

In fact, (3) is a corollary of (1)  using the main result of     \cite{[13]}, so the main work is  in proving (1)
which involves (2) as an intermediate step.

\medskip

 Note that the class of 1-dimensional Noetherian Zariski structures is essentially classifiable by the main result of \cite{[2]}, and in particular the class contains no instances of structures obtained by the proper Hrushovski construction.   The class of analytic Zariski structures, in contrast, is consistent with  Hrushovski's construction and at the same time, by the result above, has excellent model-theoretic properties. 
 This gives a hope for a classification theory based on the relevant notions.  

However, it must be mentioned that some natural questions in this context are widely open. In particular, we have no classification for presmooth
analytic Zariski groups (with the graph of multiplication analytic).  It is not known if a 1-dimensional irreducible presmooth analytic group has to be  abelian.  See related analysis of groups in  \cite{HLS}.

\medskip

{\bf Acknowledgment.}
I want to express my thanks to Assaf Hasson who saw a very early version of this work and made many useful comments, and also to Levon Haykazyan who through his own contributions to the theory of quasiminimality kept me  informed in the recent developments in the field.

\section{Analytic $L$-Zariski structures}
\lb{anzar}

Let $\M=(M, L)$ be a structure with primitives (basic relations) $L.$
 
We introduce a topology on $M^n,$ for all $n\ge 1,$ by declaring a subset $P\subs M^n$ closed if
there is an $n$-type $p$ consisting of quantifier-free positive formulas with parameters in  $M$
such that $$P=\{ a\in M^n:\ \M \vDash p(a)\}.$$ 

In other words, the sets defined by atomic $L(M)$-formulae   form a basis for the topology.

 We say $P$ is $L$-closed if  $p$ is over $\emptyset.$

\rrem Note that it follows that projections $$\prs_{i_1,\dots,i_m}:M^n\to M^m, \ \ \langle x_1, \ldots, x_n\rangle \mapsto \langle x_{i_1}, \ldots, x_{i_m}\rangle$$ 
are continuous in the sense that  the inverse image of a closed set under a projection
is closed. Indeed, $\prs\inv S=S\times M^{n-m}.$
\erem

We write $X\sbo V$ to say that $X$ is open in $V$ and
$X\sbc V$ to say  it is closed.

We say that $P\subs M^n$ is {\bf constructible}  
if $P$ is a finite union of some sets $S,$ such that $S\sbc U\sbo M^n.$

A subset $P\subs M^n$ will be called {\bf projective} if $P$ is a union of finitely many sets of the form $\pr S,$
for
some $S\sbc U\sbo M^{n+k},$ $\pr: M^{n+k}\to M^n.$

Note that any set $S$ such that $S\sbc U\sbo M^{n+k},$ is
constructible, a projection of a constructible set is projective and that any constructible set is projective.

\subsection{Dimension}
To any nonempty 
projective $S$
a non-negative integer  $\dim S,$ called
{\bf the dimension of } $S,$  is
attached.

\medskip

We assume:\\

(SI) ({\bf strong irreducibility}) for an irreducible set $S\sbc U\sbo M^n$
(that is $S$ is not a union of two proper closed subsets) and its closed subset $S'\sbc S,$ $$\dim S'=\dim S\Rightarrow S'=S;$$

(DP) ({\bf dimension of points}) for a nonempty projective  $S,$ 
$\dim S=0$ if and only if $S$ is at most countable.\\

(CU) ({\bf countable unions}) If $S= \bigcup_{i\in \N}S_i,$ all projective, then
$\dim S=\max_{i\in \N}\dim S_i;$ \\

(WP) ({\bf weak properness})
given an irreducible $S\sbc U\sbo M^n$ and $F\sbc V\sbo M^{n+k}$ with the 
projection $\pr:M^{n+k}\to M^n$ such that $\pr F\subs S$ and 
$\dim \pr F=\dim S,$ there exists $D\sbo S$
such that $D\subs \pr F.$ 
\\

\subsection{Remark}
{\rm (CU)} in the presence of
{\rm (DCC)} implies the {\em essential uncountability property} {\rm (EU)} 
usually assumed for Noetherian Zariski structures.

\medskip

We postulate further,  for an irreducible $S\sbc U\sbo M^{n+k}:$

\begin{enumerate}


\item[(AF)]  
$\dim \pr S=\dim S -   \min_{u\in \pr S} \dim (\pr\inv(u)\cap S);$

\item[(FC)]  
The set
$\{ a\in \pr S: \dim (\pr\inv(a)\cap S)\ge m\}$
 is of the form $T \cap \pr S$ for some  constructible $T,$ and 
 there exists an open set $V$
such that $V\cap\pr S\neq \emptyset$ and 
$$\min_{ a\in \pr S}
 \dim(\pr^{-1} (a)\cap S)= \dim(\pr^{-1} (v)\cap S),\mbox{ for any }v\in V\cap\pr(S).$$ 
\end{enumerate}

The following helps to understand the dimension of projective sets.
\bl \lb{dimpr}
Let $P=\pr S\subs M^n,$ for $S$ irreducible constructible, and $U\sbo M^n$
with $P\cap U\neq \emptyset.$  Then 
$$\dim P\cap U=\dim P.$$
\el
\pf We can write $P\cap U=\pr S'=P',$ where $S'=S\cap \pr\inv U$ constructible irreducible, $\dim S'=\dim S$ by (SI). By (FC), there is $V\sbo M^n$ such that for all $c\in V\cap P,$
$$\dim \pr\inv(c)\cap S=\min_{a\in P} \dim \pr\inv(a)\cap S=\dim S-\dim P.$$
Note that $\pr\inv U\cap \pr\inv V\cap S\neq \emptyset,$ since $S$ is 
irreducible. Taking $s\in \pr\inv U\cap \pr\inv V\cap S$ and $c=\pr s$ 
we get, using (AF) for $S',$
 $$\dim \pr\inv(c)\cap S'=\dim \pr\inv(c)\cap S=\min_{a\in P'} \dim \pr\inv(a)\cap S=\dim S-\dim P'.$$
So, $\dim P'=\dim P.$
\qed

\subsection{Analytic subsets}

\medskip
A subset $S,$ $S\sbc U\sbo M^n,$ is called
{\bf analytic in } $U$ if for every $a\in S$
there is an open $V_a\sbo U$ such that $a\in V_a$ and
$S\cap V_a$ is the union of finitely many closed in $V_a$ 
irreducible subsets.

\index{analytic subsets}

We postulate the following properties

\begin{enumerate}

\item[(INT)] ({\bf Intersections})  If $S_1, S_2\sba U$  are irreducible then
$S_1\cap S_2$ is analytic in $U;$\\

\item[(CMP)] ({\bf Components})  If $S\subs_{an} U$ and $a\in S,$ a closed point,  then there is
    $S_a\subs_{an} U,$ a finite union of irreducible analytic subsets
    of $U,$  and some $S'_a\subs_{an} U$ such that $a\in
    S_a\stmn S'_a$ and
    $S=S_a\cup S'_a;$\\

Each of the irreducible subsets of $S_a$ above is called an {\bf 
irreducible component of $S$
containing $a.$}

\item[(CC)] ({\bf Countability of the number of components})   Any $S\sba U$ is a 
union of at most countably many  irreducible
components.

\end{enumerate}

\subsection{Remark}
 For $S$ analytic  and $a\in \pr S,$ the fibre $S(a,M)$ 
is analytic.\label{fbr}

\bl \lb{C2} If $S\sba U$ is irreducible, $V$ open, then $S\cap V$ is
an  irreducible analytic subset of $V$ and, if non-empty, $\dim S\cap V=\dim 
S.$\el

\pf Immediate.\qed

\bl \lb{B}
\begin{itemize}
\item[(i)] \lb{B1} $\emptyset,$ any singleton and $U$ are analytic in $U;$

\item[(ii)] \lb{B3} If $S_1, S_2\sba U$  then
$S_1\cup S_2$ is analytic in $U;$

\item[(iii)] \lb{B2} If $S_1\sba U_1$ and $S_2\sba U_2,$ then
$S_1\times S_2$ is analytic in $U_1\times U_2;$

\item[(iv)] \lb{B4} If $S\sba U$ and  $V\subs U$ is open then
$S\cap V\sba V;$ 

\item[(v)] \lb{B3v} If $S_1, S_2\sba U$   then
$S_1\cap S_2$ is analytic in $U.$
\end{itemize}
\el
\pf Immediate. \qed

\df Given a subset  $S\sbc U\sbo\ M^n $ we define the notion of the {\bf 
analytic rank} of $S$ in $U,$ \ $\ark_U(S),$ which is a natural number 
satisfying \index{analytic rank $\ark$}
\begin{enumerate}
\item\lb{ark1} $\ark_U(S)=0$ iff $S=\emptyset;$

\item \lb{ark2} $\ark_U(S)\le k+1$ iff there is a set $S'\sbc\ S$ such that
$\ark_U(S')\le k$ and with the set $S^0=S\smin S'$ being analytic in $U\smin 
S'.$

\end{enumerate}

\edf

Obviously, any nonempty analytic subset of $U$ has analytic rank $1.$

\medskip

The next assumptions guarantees that the class of analytic subsets 
explicitly determines the class of closed subsets in $M.$\\

(AS) {\bf [Analytic stratification]} For any $S\sbc U\sbo M^n,$ $\ark_U S$ is
defined and is finite.

\medskip

We will justify this non obvious property later in  \ref{gan1} and \ref{gan2}.

\bl \lb{dimcl} For any $S\sbc U\sbo M^n,$
$$\dim\pr S+ \min_{a\in \pr S}\dim \pr\inv(a)\cap S\ge\dim S.$$
\el
\pf We use (AS) and prove the statement by induction on $\ark_US\ge 1.$

For $\ark_U S=1,$ $S$ is analytic in $U$ and so by (CC) is the union of
countably many irreducibles $S_i.$ By (AF)
$$\dim\pr S_i+ \min_{a\in \pr S_i}\dim \pr\inv(a)\cap S_i\ge\dim S_i$$
and so by (CU) lemma follows.\qed

\subsection{Presmoothness}
The following  property (which we are not going to use in the context of the present paper) is relevant.

\medskip

(PS) {\bf [Presmoothness]} If $S_1, S_2\sba U\sbo M^n$ and $S_1,S_2$ and $U$
irreducible, then for any irreducible
component $S_0$ of  $S_1\cap S_2$ $$\dim S_0\ge \dim S_1+\dim S_2
-\dim U.$$

\df
An $L$-structure $\M$ is said to be  {\bf analytic $L$-Zariski}  if
\begin{itemize}
\item $\M$
  satisfies 
(SI), (WP), (CU), (INT), (CMP),(CC), (AS);
\item
the expansion  $\M^\sharp$ of $\M$ to the language $L(M)$ (names for points in $M$ added) satisfies all the above  with the dimension extending the one for $\M;$
\item $\M^\sharp$ also satisfies  
(AF) and (FC)   with $V$ in (FC) being $L$-definable whenever  $S$ is.
\end{itemize}
An analytic Zariski structure will be called {\bf presmooth} if it has the
presmoothness property
(PS).\edf

\section{Model theory of analytic
Zariski structures}\lb{qe}

\bl Let $\M$ be analytic $L$-Zariski and assume $L$ is countable.
There is a countable $\M_0\preccurlyeq \M$ 
such that for   
any $L$-closed set $S$ any   irreducible component $P$ of $S$ is $L$-definable ($L$-closed); 
\el
\pf Use the standard L\"owenheim - Skolem downward arguments. \qed

We call $M_0$ a core substructure (subset) of $\M.$

\subsection{Assumption} By extending $L$ to $L(M_0)$ we  assume that   the set of $L$-closed points is the core subset.

\df For finite $X\subs M$  we define the {\bf predimension}
\be \delta(X)=\min\{ \dim S:\ \vec{X}\in S,\ S\sba U\sbo M^n,\ S \mbox{ is 
$L$-constructible}  \},\ee
{\bf relative predimension} for finite $X,Y\subs M$ 
\be \label{dlt} \delta(X/Y)=\min\{ \dim S:\ \vec{X}\in S,\ S\sba U\sbo M^{n},\ S \mbox{ is 
$L(Y)$-constructible}  \},\ee

and {\bf dimension of $X$}
$$\dd(X)=\min\{ \delta(XY):\ \mbox{ finite } Y\subset M\}.$$

We call a minimal $S$ as in (\ref{dlt}) an {\bf analytic locus of $X$ over $Y.$}

\medskip  

For $X\subs {M}$ finite, we say that $X$ is {\bf
self-sufficient}
and write $X\le {M},$ if $\dd(X)=\delta(X).$ 

For infinite $A\subs M$ we say $A\le M$ if for any finite $X\subs A$ there is
a finite $X\subs X'\subs A$ such that $X'\le M.$
\edf

\subsection{Assumption} $\dim M=1$ and $M$ is irreducible.

Note that we then have
$$0\le\delta(Xy)\le \delta(X)+1,\mbox{ for any }y\in M,$$
since $\vec{Xy}\in S\times M.$

\bl\lb{prk} Given $F\sba U\sbo M^k,$  $\dim F>0,$ there is $i\le k$
such that for $\prs_i: (x_1,\ldots,x_k)\mapsto x_i,$ 
$$\dim\prs_iF>0.$$
\el
\pf Use (AF) and induction on $k.$\qed

\bp \lb{qe1} Let $P=\pr S,$ for some $L$-constructible
$S\sba U\sbo M^{n+k},$ $\pr: M^{n+k}\to M^n.$  Then
\be \lb{form}\dim P=\max \{ \dd(x):\ x\in P(M)\}.\ee
Moreover, this formula is true when $S\sbc U\sbo M^{n+k}.$

\ep
\pf We use induction on $\dim S.$ 

We first note that by induction on $\ark_US,$
if (\ref{form}) holds for all analytic $S$ of
dimension less or equal to $k$ then it holds for all closed  $S$ of
dimension less or equal to $k.$

The statement is obvious for $\dim S=0$ and so we assume that $\dim
S>0$ 
and for all analytic $S'$ of lower dimension the statement is true.

By (CU) and (CMP) we may
assume that $S$ is irreducible. Then by (AF) 
\be \lb{Sc}\dim P=\dim S-\dim S(c,M)\ee for any $c\in P(M)\cap V(M)$ (such that $S(c,M)$ is of minimal 
dimension)
for some open $L$-constructible $V.$ 

Claim 1. It suffices to prove the statement of the proposition for the 
projective set
$P\cap V',$ for some $L$-open $V'\sbo M^n.$ 

Indeed, 
$$P\cap V'=\prs(S\cap \pr\inv V'),\ \ S\cap \pr\inv V'\sbc \pr\inv V'\cap U\sbo 
M^{n+k}.$$
And $P\smin V'=\prs(S\cap T)$, $T=\pr\inv(M^n\smin V')\in L.$ So,
$P\smin V'$ is the projection of a proper analytic subset, of lower dimension.
By induction, for $x\in P\smin V',$ $\dd(x)\le \dim P\smin V'\le \dim P$ and hence, using \ref{dimpr},
$$\dim P\cap V'=\max \{ \dd(x):\ x\in P\cap V'\}\Rightarrow
\dim P=\max \{ \dd(x):\ x\in P\}.$$      

Claim 2. The statement of the proposition holds if  $\dim S(c,M)=0$ in (\ref{Sc}). 
 
Proof. Given $x\in P$ choose a tuple $y\in M^k$ such that
    $S(x\concat y)$ holds. Then $\delta(x\concat y)\le \dim S.$ So we
have
$\dd(x)\le \delta(x\concat y)\le \dim S=\dim P.$

It remains to notice that there exists $x\in P$ such that
$\dd(x)\ge\dim P.$

Consider the   $L$-type
$$x\in P\ \& \{ x\notin R: \ \dim R\cap P<\dim P\mbox{ and $R$ is
projective}\}.$$
This is realised in $M,$ since
otherwise $P=\bigcup_{R} (P\cap R)$ which would contradict
(CU) because $(M_0,C_0)$ is countable.

For such an $x$ let $y$
be a tuple in $M$ such that $\delta(x\concat y)=\dd(x).$ By definition there
exist $S'\sba U'\sbo M^m$ such that $\dim S'=\delta(x\concat y).$ Let
$P'=\pr S',$
the projection into $M^n.$ By our choice of $x,$ $\dim P'\ge \dim P.$
But $\dim S'\ge \dim P'.$ Hence, $\dd(x)\ge \dim P.$ Claim proved.\\
 
Claim 3. There is a $L$-constructible $R\sba S$ such that all the
fibres $R(c,M)$ of the projection map $ R\to \pr R$ are $0$-dimensional and
$\dim \pr R=\dim P.$

Proof.
We have by construction $S(c,M)\subs M^k.$ 
Assuming $\dim S(c,M)>0$ on every open subset
we  show that there is a
$b\in M_0$ such that (up to the order of coordinates) 
$\dim S(c,M)\cap \{ b\}\times M^{k-1}<\dim S(c,M),$ 
 for all $c\in P\cap V'\neq \emptyset,$ for some
open $V'\subs V$ and $\dim \pr S(c,M)\cap \{ b\}\times M^{k-1}=\dim
P.$ 
By induction on $\dim S$ this will prove the claim. 

To find such a $b$ choose $a\in P\cap V$ and note that by \ref{prk},
up to the order of coordinates, $\dim\prs_1S(a,M)>0,$
where $\prs_1: M^k\to M$ is the projection
on the first coordinate.

Consider the projection $\prs_{M^n,1}: M^{n+k}\to M^{n+1}$ and
the set $\prs_{M^n,1}S.$
By (AF) we  have
$$\dim  \prs_{M^n,1}S= \dim P+\dim \prs_1S(a,M)=\dim P+1.$$

Using (AF) again for the projection $\prs^1:M^{n+1}\to M$ with the
fibres $M^n\times \{ b\}, $
  we get, for all $b$ in some open subset of $M,$
$$1\ge \dim  \prs^1\prs_{M^n,1}S= \dim  \prs_{M^n,1}S-\dim[\prs_{M^n,1}S]\cap
[M^n\times \{ b\}]=$$ $$=\dim P+1-\dim[\prs_{M^n,1}S]\cap
[M^n\times \{ b\}].$$
Hence $\dim[\prs_{M^n,1}S]\cap
[M^n\times \{ b\}]\ge\dim P,$ for all such $b,$ which means that the 
projection of
the set $S_b=S\cap (M^n\times \{ b\}\times M^{k-1})$ on $M^n$ is of dimension
$\dim P,$ which finishes the proof if $b\in M_0.$ 
But $\dim S_b=\dim S-1$ for all $b\in M\cap V',$ some $L$-open
$V',$ so 
for any $b\in M_0\cap V'.$ 
The latter is not empty since $(M_0,L)$ is a core substructure.
This proves the claim.\\

Claim 4. Given $R$ satisfying Claim 3, 
$$P\smin \pr R\subs \pr S',\mbox{ for some } S'\sbc S,\ \dim S'<\dim S.$$
Proof. Consider the cartesian power $$M^{n+2k}=\{ x\concat y\concat z:
x\in M^n,\ y\in M^k,\ z\in M^k\}$$
and its $L$-constructible subset $$R\& S:=\{ x\concat y\concat z:\ \ x\concat z\in R\
\&\  x\concat y\in S\}.$$
Clearly $R\& S\sba W\sbo M^{n+2k},$ for an appropriate
$L$-constructible $W.$ 

Now notice that the fibres of the projection $\prs_{xy}: x\concat y\concat z\mapsto
x\concat y$ over $\prs_{xy}R\& S$ are $0$-dimensional and so, for some
irreducible component $(R\& S)^0$ of the analytic set $R\& S,$ 
$\dim\prs_{xy}(R\& S)^0=\dim S.$ Since $\prs_{xy}R\& S\subs S$ and 
$S$ irreducible, we get by (WP)
$D\subs \prs_{xy}R\& S$ 
for some $D\sbo S.$ Clearly $$\pr R=\pr\prs_{xy}R\& S\sups \pr D$$
and $S'=S\smin D$ satisfies the requirement of the claim.\\

Now we complete the proof of the proposition: By Claims 2 and 3 
$$\dim P=\max_{x\in \pr R} \dd(x).$$
By induction on $\dim S,$ using Claim 4, for all $x\in P\smin \pr R,$
$$\dd(x)\le \dim \pr S'\le \dim P.$$
The statement of the proposition follows. 
\qed

In what follows a $L$-substructure of $\M$ is a $L$-structure on a subset $N\sups M_0.$ Recall that $L$ is purely relational.
 
Recall the following well-known fact, see \cite{[7*]}.

\subsection{Karp's characterisation of $\equiv_{\infty,\omega}$} \label{Fact} 
Given $a,a'\in M^n$ 
the  $L_{\infty,\omega}(L)$-types of the two $n$-tuples in $\M$ are equal
if and only if 
they are {\em back and forth equivalent} that is  there is a
nonempty set $I$ of isomorphisms of $L$-substructures  of $\M$ such
that $a\in \dom f_0$ and $a'\in {\rm Range}\, f_0,$ for some $f_0\in I,$ and

(forth) for every $f\in  I$ and $b\in M$ there is a $g\in I$ such that $f\subs  g$ and $b\in  \dom g;$

(back) For every $f\in I$ and $b'\in M$ there is a $g\in I$ such that $f\subs  g$
 and $b'\in {\rm Range}\, g.$ 

\df For $a\in M^n,$
the {\bf projective type of $a$
over $M$} is $$\{ P(x):\ a\in  P,\ P \mbox{ is a projective set over
$L$}\}\cup$$
$$\cup \{ \neg P(x):\ a\notin  P,\ P \mbox{ is a projective set over
$L$}\}.$$

\edf

\bl \lb{XX'} Suppose  $X\le {M},$  $X'\le {M}$ and the (first-order)
quantifier-free $L$-type of $X$
is equal to that of $X'.$  Then the  $L_{\infty,\omega}(L)$-types of 
$X$ and $X'$ are equal.
\el

\pf We are going to construct a back-and-forth system for $X$ and $X'.$

Let $S_X\sba V\sbo M^n,$ $S_X$ irreducible, all $L$-constructible, and such 
that $X\in S_X(M)$ and
$\dim S_X=\delta(X).$

Claim 1. The quantifier-free
$L$-type of $X$ (and $X'$) is determined by formulas equivalent to 
$S_X\cap V',$ for
$V'$ open such that $X\in V'(M).$

Proof. Use  the stratification of closed sets (AS) to choose $L$-constructible  $S\sbc U\sbo M^n$
such that $X\in S$ and
$\ark_US$ is minimal. Obviously then $\ark_US=0,$ that is $S\sba U\sbo M^n.$
Now $S$ can be decomposed into irreducible components, so we may choose $S$
to be irreducible. Among all such $S$ choose one which is of minimal
possible dimension. Obviously $\dim S=\dim S_X,$ that is we may assume that
$S=S_X.$
Now clearly any constructible set $S'\sbc U'\sbo M^n$ containing $X$
must satisfy $\dim S'\cap S_X\ge \dim S_X,$ and this condition is also sufficient for $X\in S'.$  
\\

Let $y$ be an element of $M.$ We want to find a finite $Y$
containing $y$ and
an $Y'$ such that the quantifier-free type of $XY$ is equal to that
of $X'Y'$ and both are self-sufficient in $M.$ This, of course,
extends the partial isomorphism $X\to X'$ to $XY\to X'Y'$ and will prove
the lemma.

We choose $Y$ to be a minimal set containing $y$ and such that $\delta(XY)$
is also minimal, that is $$1+\delta(X)\ge\delta(Xy)\ge\delta(XY)=\dd(XY)$$ and 
$XY\le {M}.$

We have two cases: $\delta(XY)=\dd(X)+1$ and $\delta(XY)=\dd(X).$ In
the first case $Y=\{ y\}.$ By Claim 1 the quantifier-free $L$-type 
$r_{Xy}$ of
$Xy$ is determined by the formulas of the form
$(S_X\times M)\smin T,$\
$T\subs_{cl} M^{n+1},$ $T\in L,$ $\dim T<\dim (S_X\times M).$

Consider $$r_{Xy}(X',M)=\{ z\in M: X'z\in (S_X\times M)\smin T,\ \dim
T<\dim S_X, \mbox{ all } T
\}.$$

We claim that $r_{Xy}(X',M)\neq \emptyset.$
Indeed, otherwise $M$ is the union of countably many sets of
the form $T(X',M).$ But
the fibres $T(X',M)$ of $T$ are of dimension $0$ (since otherwise
$\dim T=\dim S_X+1,$ contradicting the definition of the $T$). This is
impossible, by (CU).

Now we choose $y'\in r_{Xy}(X',M)$ and this is as required.

In the second case,
by definition, there is an irreducible $R\sba U\sbo M^{n+k},$ $n=|X|, 
k=|Y|,$
such that $XY\in R(M)$ and $\dim R=\delta(XY)=\dd(X).$ We may assume
$U\subs V\times M^k.$

Let $P=\pr R,$ the
projection into $M^n.$ Then $\dim P\le \dim R.$ But also $\dim P\ge
\dd(X),$
by \ref{qe1}.  Hence, $\dim R=\dim P.$ On the other hand, $P\subs S_X$
and $\dim S_X=\delta(X)=\dim P.$ By axiom (WP) we have 
$S_X\cap V'\subs P$ for some $L$-constructible open $V'.$

Hence  $X'\in S_X\cap V' \subs P(M),$ for
$P$ the projection of an irreducible
analytic set $R$ in  the $L$-type of $XY.$ By Claim 1 the
quantifier-free $L$-type of $XY$ is of the form
$$r_{XY}=\{ R\smin T: T\sbc R,\ \dim T<\dim R\}.$$
Consider $$r_{XY}(X',M)=\{ Z\in M^k: X'Z\in R\smin T,\ T\sbc R,\ \dim T<\dim 
R\}.$$

We claim again that $r_{XY}(X',M)\neq \emptyset.$
Otherwise the set $R(X',M)=\{ X'Z:\ R(X'Z)\}$
is the union of countably many subsets of the form $T(X',M).$
But $\dim T(X',M)<\dim R(X',M)$ as above, by (AF).

Again, an $Y'\in  r_{XY}(X',M)$ is as required.\qed

\bc \lb{cbly}There is at most countably many  $L_{\infty,\omega}(L)$-types of tuples
$X\le M.$
\ec

Indeed, any such type is determined uniquely by the choice of a $L$-constructible
 $S_X\sba U\sbo M^n$ such that $\dim S_X=\dd(X).$ 

\bl\lb{6.5.10} Suppose, for finite $X,X'\subs {M},$ the projective $L$-types of $X$
and $X'$ coincide. Then the  $L_{\infty,\omega}(L)$-types of the tuples are equal.
\el

\pf Choose finite $Y$ such that $\dd(X)=\delta(XY).$ Then $XY\le M.$
Let $XY\in S\sba U\sbo M^n$ be $L$-constructible and
such that $\dim S$ is minimal possible, 
that is $\dim S=\delta(XY).$ We may assume that $S$ is irreducible.
Notice that for every proper closed $L$-constructible $T\sbc U,$
$XY\notin T$ by dimension considerations.

By assumptions of the lemma $X'Y'\in S,$ for some $Y'$ in $M.$
We also have $X'Y'\notin T,$ for any $T$ as above, since otherwise
a projective formula would imply that $XY''\in T$ for some $Y'',$ 
contradicting that $\dd(X)>\dim T.$ 

We also have $\delta(X'Y')=\dim S.$ But for no finite $Z'$ it is
possible that $\delta(X'Z')<\dim S,$ for then again a projective
formula will imply that $\delta(XZ)<\dim S,$ for some $Z.$ 

It follows that $X'Y'\le M$ and 
the quantifier-free types of $XY$ and $X'Y'$ coincide, 
hence the  $L_{\infty,\omega}(L)$-types are equal, by \ref{XX'}.\qed   

\df Set, for finite $X\subs M,$
$$\cl_{L}(X)=\{ y\in M:\ \dd(Xy)=\dd(X)\}.$$

\edf

We fix $L$ and omit the subscript below.

\bl \lb{ccbl}   The following two conditions are equivalent
\begin{itemize}
\item[(a)]$b\in \cl(A),$ for $\vec{A}\in M^n;$ 
\item[(b)]  $b\in P(\vec{A},M)$ for some
projective first-order $P\subs M^{n+1}$ such that  $P(\vec{A},M)$ is at most countable. 
\end{itemize}

In particular, $\cl(A)$ is countable for any finite $A.$
\el

\pf Let $d=\dd(A)=\delta(AV),$ and $\delta(AV)$ is minimal for all possible finite
$V\subs M.$ So by definition $d=\dim S_0,$ some  analytic irreducible $S_0$ such that
$\vec{AV}\in S_0$ and $S_0$ of minimal dimension. This corresponds to a 
$L$-definable relation
$S_0(x,v),$ where $x,v$  strings of variables of length $n,m$ 

First assume (b), thais is that $b$ belongs to a countable $P(\vec{A},M).$ By definition
$$P(x,y)\equiv \exists w\, S(x,y,w),$$
for some analytic $S\subs M^{n+1+k},$ some tuples $x,y,w$  of variables of length $n,1$ and 
$k$ respectively, and the fibre $S(\vec{A},b,M^k)$ is nonempty. We also assume that $P$ and $S$ are 
of minimal dimension, answering this description.
By (FC), (AS) and minimality we may choose
$S$ so that $\dim S(\vec{A},b,M^k)$ is minimal among all the fibres
$S(\vec{A'},b',M^k).$

Consider the analytic set  $S^\sharp\subs M^{n+m+1+k}$ given by 
$S_0(x,v)\,\&\, S(x,y,w).$ By (AF), considering the projection of the set on $(x,v)$-coordinates, 
$$\dim S^\sharp\le \dim S_0+\dim S(\vec{A},M,M^k),$$
since $S(\vec{A},M,M^k)$ is a fibre of the projection. Now we note that
by countability $\dim S(\vec{A},M,M^k)=\dim S(\vec{A},b,M^k),$ so
 $$\dim S^\sharp\le \dim S_0+\dim S(\vec{A},b,M^k).$$
Now the projection  $\prs_w S^\sharp$ along $w$ (corresponding to $\exists w\, S^\sharp$)
has fibres of the form  $S(\vec{X},y,M^k),$ so by (AF)
$$\dim \prs_w S^\sharp\le \dim S_0= d.$$
Projecting further along $v$ we get
$\dim \prs_v\prs_w S^\sharp\le d,$ but $\vec{A}b\in \prs_v\prs_w S^\sharp$
so by Proposition~\ref{qe1} $\dd( \vec{A}b)\le d.$ The inverse inequality holds by 
definition, so the equality holds. This proves that $b\in \cl(A).$

Now assume (a), that is $b\in \cl(A).$ So, $\dd(\vec{A}b)=\dd(\vec{A})=d.$
By definition there is a projective set $P$ containing $\vec{A}b,$ defined by 
the formula $\exists w\, S(x,y,w)$ for some analytic $S,$ $\dim S=d.$  
Now $\vec{A}$ belongs to the projective set $\prs_yP$  (defined by the formula
$\exists y\exists w\,S(x,y,w)$) so by Proposition~\ref{qe1}
$d\le \dim \prs_y P,$ but $\dim \prs_yP\le \dim P\le \dim S=d.$
Hence all the dimensions are equal and so, the dimension of the generic fibre is $0$. 
We may assume, as above, without loss of generality that all fibres are of minimal dimension, so 
$$\dim S(\vec{A},M, M^k)=0.$$
Hence, $b$ belongs to a 0-dimensional set $\exists w\, S(\vec{A},y,w),$ which is projective and countable. \qed

\bl \label{QM1} Suppose $b\in \cl(A)$ and the projective type of $\vec{A}b$ is equal to that of  $\vec{A'}b'.$ Then $b'\in \cl(A').$  
\el

\pf First note that, by (FC) and (AS), for analytic $R(u,v)$ and its fibre $R(a,v)$ of minimal dimension one has 
$$\mathrm{tp}(a)=\mathrm{tp}(a') \Rightarrow
\dim R(a,v)=  \dim R(a',v).$$

By the second part of the proof of \ref{ccbl} the assumption of the lemma implies that for some analytic $S$ we have $\vDash\exists w  S(\vec{A},b, w)$ and
 $\dim S(\vec{A},M, M^k)=0.$ Hence  $\vDash\exists w  S(\vec{A'},b', w)$ and
 $\dim S(\vec{A'},M, M^k)=0.$ But this immediately implies $b'\in \cl(A').$ \qed

\bl \lb{preg}

 (i)
$$\cl(\emptyset)=\cl(M_0)=M_0.$$ 
(ii) Given finite $X\subs M,$ $y,z\in M,$ 
$$z\in \cl(X,y)\smin \cl(X)\Rightarrow y\in \cl(X,z).$$
(iii) $$\cl(\cl(X))=\cl(X).$$
\el
\pf (i) Clearly $M_0\subs \cl(\emptyset),$ by definition.

We need to show the converse, that is if $\dd(y)= 0,$ for $y\in M,$ then $y\in M_0.$
By definition $\dd(y)=\dd(\emptyset)=\min\{\delta(Y):y\in Y\subset M\}= 0.$
 So, $y\in Y,$ $\vec{Y}\in S\sba U\sbo M^n,$ $\dim S=0.$ The irreducible components of $S$ are closed points (singletons) and $\{ \vec{Y}\}$ is one of them, so must be in $M_0,$ hence   $y\in M_0.$

(ii) Assuming the left-hand side of (ii), 
$\dd(Xyz)=\dd(Xy)>\dd(X)$ and  
 $\dd(Xz)>\dd(X).$  
 By the definition of $\dd$ then, 
$$\dd(Xy)=\dd(X)+1=\dd(Xz),$$
so  $\dd(Xzy)=\dd(Xz),$ $y\in \cl(Xz).$

(iii) Immediate by \ref{ccbl}.
\qed

Below, if not stated otherwise, we use the language $L^\exists$ the primitives of which correspond to relations $\exists$-definable in $\M.$
Also, we call a {\bf submodel of $\M$} any $L^\exists$-substructure  closed under $\cl.$

\bt \lb{proj} (i) Every  $L_{\infty,\omega}(L)$-type realised in $\M$ is equivalent to a projective
type, that is a type consisting of existential (first-order) formulas and the
negations of existential formulas.

(ii) There are only countably many $L_{\infty,\omega}(L)$-types realised in
$\M.$

(iii)  
$(M,L^{\exists})$ is quasi minimal $\omega$-homogeneous over countable submodels, that is the following hold:
 \begin{itemize}
 
\item[(a)]
  for any countable (or empty) submodel $G$ and any $n$-tuples $X$ and $X',$ both $\cl$-independent over $G,$ 
a bijection $\phi:X\to X'$ is a $G$-monomorphism;
\item[(b)]
 given any $G$-monomorphism 
$\phi: Y\to Y'$ for finite tuples $Y,Y'$ in $M$ and given
 a $z\in M$
we can extend $\phi$ so that $z\in \dom \phi.$
 \end{itemize}
\et

\pf (i) Immediate from \ref{6.5.10}. 

(ii) By
\ref{cbly} there are only countably many types of finite tuples $Z\le M.$
Let $N\subs M_0$ be a countable subset of $M$ such that any finite
$Z\le M$ is $L_{\infty,\omega}(L)$-equivalent to some tuple in $N.$
Every finite tuple $X\subset M$ can be extended to $XY\le M,$ so there is
a $L_{\infty,\omega}(L)$-monomorphism $XY\to N.$ This monomorphism
identifies the $L_{\infty,\omega}(L)$-type of $X$ with one of a tuple in $N,$ hence there are no more than countably many such types.

(iii) Lemma~\ref{preg} proves that $\cl$ defines a pregeometry on $M.$

Consider first (a). Note that $GX\le M$ and $GX'\le M$ and so the types of
$X$ and $X'$ over $G$ are $L$-quantifier-free. But there is no proper 
$L$-closed subset $S\sbc M^n$ such that $\vec{X}\in S$ or  $\vec{X'}\in S.$
Hence the types are equal.

For (b) just use the fact that the $G$-monomorphism by our definition
preserves $\exists$-formulas, so by \ref{6.5.10}
complete $L_{\infty,\omega}(L(G))$-types of
$X$ and $X'$  coincide, so by \ref{Fact} $\phi$ can be extended.  
\qed

\medskip

\bt \label{main} $M$ is a quasiminimal pregeometry structure (see \cite{[13]}). In other words, the following properties of $\M$ hold:

(QM1) The pregeometry $\cl$ is determined by the language. That is, if $\mathrm{tp}(x,Y ) = \mathrm{tp}(x',Y')$ , then
$x\in \cl(Y)$ if and only $x'\in \cl(Y').$ (Here the types are first order).

(QM2) The structure $\M$ is infinite-dimensional with respect to $\cl.$

(QM3) (Countable closure property). If $X\subset M$ is finite, then $\cl(X)$ is countable.

(QM4) (Uniqueness of the generic type). Suppose that $H,H'\subset M$ are countable closed subsets,
enumerated such that $\mathrm{tp}(H)=\mathrm{tp}(H').$
 If $y\in M \setminus H$ and $y'\in M\setminus  H',$
 then $\mathrm{tp}(H,y)=\mathrm{tp}(H',y').$

(QM5) ($\omega$-homogeneity over closed sets and the empty set). Let$H,H'\subset M$ be countable
closed subsets or empty, enumerated such that $\mathrm{tp}(H)=\mathrm{tp}(H'),$ and let $Y, Y'$ be finite tuples
from $M$ such that $\mathrm{tp}(H,Y)=\mathrm{tp}(H',Y'),$
 and let $z\in \cl(H,Y).$ Then there is $z'\in  M$ such that
$\mathrm{tp}(H,Y,z)=\mathrm{tp}(H',Y',z').$

\et

\pf $(M,\cl)$ is a pregeometry by \ref{preg}.
 (QM1) is proved in \ref{QM1}. (QM3) is \ref{ccbl} and (QM2) follows from (QM3)  and (CU). (QM4)\&(QM5) is \ref{proj}(iii). 
 \qed
 
 \medskip 

Now we  define an {\em abstract elementary class} $\A$ 
associated with $\M.$ We follow \cite{Znemt},   for this construction. Similar construction was used in \cite{[13]}.
 
Set  
$$\A_0(\M)=\{ \mbox{ countable $L^{\exists}$-structures } {\rm N}:\ {\rm N}\cong {\rm N}'\subs \M,\ \cl(N')=N'\}$$
and define embedding ${\rm N}_1\preccurlyeq {\rm N}_2$ in the class as an 
$L^{\exists}$-embedding $f:{\rm N}_1\to {\rm N}_2$ such
that there are isomorphisms $g_i:{\rm N}_i\to {\rm N}'_i,$ ${\rm N}'_1\subs {\rm N}'_2\subs \M,$ all embeddings commuting and $\cl(N'_i)=N'_i.$ 

Now define $\A(\M)$ to be the class of $L$-structures ${\rm H}$ with $\cl_{L}$ defined with 
respect to ${\rm H}$ and satisfying:

(i) $\A_0({\rm H})\subs \A_0(\M)$  as classes with embeddings

and 

(ii) for every finite $X\subs H$ there is ${\rm N}\in \A_0({\rm H}),$
such that $X\subs N.$ \\ 

Given ${\rm H}_1\subs {\rm H}_2,$ ${\rm H}_1,{\rm H}_2\in \A(\M),$ we define ${\rm H}_1\preccurlyeq {\rm H}_2$ to hold in the class, if
for every finite $X\subs H_1,$ $\cl(X)$ is the same in ${\rm H}_1$
and ${\rm H}_2.$ 
More generally, for ${\rm H}_1,{\rm H}_2\in \A(\M)$ we define
${\rm H}_1\preccurlyeq_f {\rm H}_2$ to be an embedding $f$ such that there are isomorphisms
${\rm H}_1\cong {\rm H}'_1,$ ${\rm H}_2\cong {\rm H}'_2$ such that
${\rm H}'_1\subs {\rm H}'_2,$ all embeddings commute,
 and ${\rm H}_1\preccurlyeq {\rm H}_2.$

\bl $\A(\M)$ is closed under the unions of ascending $\preccurlyeq$-chains. 
\el

\pf Immediate from the fact that for infinite $Y\subs M,$
$$\cl(Y)=\bigcup\{ \cl(X):\ X\subs_{\rm finite} Y\}.$$\qed

\bt \lb{aec} 
 The class $\A(\M)$ contains structures of any infinite cardinality
and is categorical in uncountable cardinals.
\et

\pf This follows from \ref{main} by the main result of \cite{[13]}.

\qed\\

\bp  
 Any uncountable ${\rm H}\in \A(\M)$ is an analytic Zariski structure in the language $L$ with parameters in ${\rm H}.$ Also ${\rm H}$ is presmooth if $\M$ is.
\ep

\pf We define $\mathcal{C}({\rm H})$ to consist of the subsets of $H^n$ of the form $P(a,H)$ for $P\in L$ of arity $k+n,$ $a\in H^k.$ The assumption (L) is obviously satisfied. 

Now note that the constructible and projective sets in $\mathcal{C}({\rm H})$ are also of the form $P(a,H)$ for some $L$-constructible or $L$-projective $P.$

Define $\dim P(a,H)=d$ if $\dim P(b,M)=d$ for some $b\in M^k$ such that
The $L^\exists$-quantifier-free types of $a$ and $b$ are equal. 
This is well-defined by (FC) and the fact that the any  $L^\exists$-quantifier-free type realised in ${\rm H}$ is also realised in $\M.$
Moreover, we have the following.

\noindent Claim. The set of $L^\exists$-quantifier-free types realised in ${\rm H}$ is equal to that realised in $\M.$

Indeed, this is immediate from the definition of the class $\A(\M),$ stability of $\A(\M)$ and the fact that the class is categorical in uncountable cardinalities.

The definition of dimension immediately implies (DP), (CU),(AF) and (FC) for ${\rm H}.$

(SI): if $P_1(a_1,H)\sbc P_0(a_0,H),$  $\dim P_1(a_1,H)=\dim P_0(a_0,H)$ and the two sets are not equal, then the same holds for  $P_1(b_1,M)$ and $P_0(b_0,M)$ for equivalent $b_0,b_1$ in $\M.$ Then, $P_0(b_0,M)$ is reducible, that is for some proper $P_2(b_2,M)\sbc P_0(b_0,M),$ $P_0(b_0,M)=P_1(b_1,M)\cup P_2(b_2,M).$
Now, by homogeneity we can choose $a_2$ in ${\rm H}$ such that
$P_0(a_0,H)=P_1(a_1,H)\cup P_2(a_2,H),$ a reducible representation.

This also shows that the notion of irreducibility is preserved by equivalent substitution of parameters. Then the same is true for the notion of analytic subset. Hence (INT), (CMP),(CC) and (PS) follow. For the same reason (AS) holds.
Next we notice that the axioms (WP) follows by the homogeneity argument. \qed

\section{Some examples} \label{s2}

\subsection{Universal covers of semiabelian varieties}

Let $\Ac$ be a semiabelian variety of dimension $d,$ e.g. $d=1$ and $\Ac$ the algebraic torus $\C^\times.$ Let $V$ be the universal cover of $\Ac,$ which classically can be identified as  a complex manifold $\C^d.$

We define a structure with a (formal)  topology on $V$
and show that this is analytic Zariski.  

By definition of universal cover there is a covering holomorphic map
 $$\exp: V\to \Ac$$
 (a generalisation of the usual $\exp$ on $\C$).

We will assume that $\A$ has no proper semiabelian subvarieties (is simple) and no complex multiplication. 

We consider the two sorted structure $(V,\Ac)$ in the language 
that has all Zariski closed subsets of $\Ac^n,$ all $n,$ the addition $+$ on $V$ and the map $\exp$ as the primitives.

This case was first looked at model-theoretically in  \cite{[11]} and the special case $\Ac=\C^\times$ in \cite{[9]},
\cite{[10]} and in the DPhil thesis \cite{[12]} of Lucy Smith.

Our aim here is to show that the structure on the sort $V$ with a naturally given formal topology is analytic Zariski.
 
\medskip

The positive quantifier-free definable subsets of $V^n,$ $n=1,2,\ldots$ form a base of a topology which we call 
 the PQF-topology. In other words

\df \lb{Lemma3.2.1} A PQF-closed subset of $V^n$ is defined as a finite union of sets of the form 
                   \be \lb{PQF1} L\cap m\cdot\ln W\ee
where $W\subs \Ac^n,$ an algebraic subvariety, and $L$ is a  $\Q$-linear subspace of $V^n,$ that is defined by a system of  equations of the form $m_1x_1+\ldots+m_nx_n=a,$ $m_i\in \Z,\  a\in V^n.$

The relations on $V$ which correspond to PQF$_\omega$-closed sets are the primitives of our language $L.$

PQF-closed subsets form a base for a topology on the cartesian powers of $V$ which will underlie the analytic Zariski structure on $V.$ 
\edf

{\bf Remark.} Among closed sets of the topology we have
 sets of the form
$$\cup_{a\in I}(S+a)$$
where $S$ is of the form (\ref{PQF1})  and  $I$ a subset of $(\ker \exp)^n.$

Slightly rephrasing the quantifier-elimination statement proved in \cite{[11]} 
Corollary 2 of section 3, we have the following result.
\bp \lb{qelim} (i) Projection of a PQF-closed set is PQF-constructible, that is
a boolean combination of  PQF-closed sets.

(ii) The image of a constructible set under exponentiation is a Zariski-constructible (algebraic) subset of $\Ac^n.$ The image of the set
 of the form (\ref{PQF1}) is Zariski closed. 
\ep


We assign {\bf dimension} to a closed set of the form (\ref{PQF1})
$$\dim L\cap m\cdot\ln W:=\dim \exp \left( L\cap m\cdot\ln W\right).$$
using the fact that the object on the right hand side is an algebraic
variety. We extend this to an arbitrary closed set assuming (CU), that is
that the dimension of a countable union is the maximum dimension of its members.This immediately gives (DP). Using \ref{qelim} we also get (WP).

The analysis of irreducibility below is more involved. Since $\exp L$ is definably and topologically isomorphic to $\Ac^k,$ 
some $k\ge 1,$ we can always reduce the analysis of a closed set of the form (\ref{PQF1}) to a one of the form $\ln W$ with $W\subs \Ac^k$ not contained in a coset of a proper algebraic subgroup. 

For such  a $W$ consider  its $m$-th ``root'' 
$$W^{\frac{1}{m}}=
\{ \la x_1,\ldots,x_n\ra\in \Ac^k: \la x^m_1,\ldots,x^m_k\ra\in W\}.$$
 Let $d=d_W(m)$ be the number of irreducible components of $ W^{\frac{1}{m}}.$ 

It is easy to see that if $d>1,$ irreducible components $W^{\frac{1}{m}}_i,$ $i=1,\ldots, d,$ of $W^{\frac{1}{m}}$ are shifts of each other by $m$-th roots of unity, and  $m\cdot \ln W^{\frac{1}{m}}_i$ are proper closed subsets of $\ln W$ of the same dimension. It follows that {\em $\ln W$ is irreducible (in the sense of {\rm (SI)}) if and only if $d_W(m)=1$ for all $m\ge 1.$ } In \cite{BGH} $W$ satisfying this condition is called {\bf Kummer generic}. If $W\subset \exp L$ for some  $\Q$-linear subspace $L\subset V^n,$ then one uses the relative version of Kummer genericity.

We say that the sequence $W^{\frac{1}{m}},$ $m\in \N,$  stops branching if the sequence $d_W(m)$ is eventually constant, that is if  $W^{\frac{1}{m}}$ is Kummer generic for some $m\ge 1$.

The following is proved for $\Ac=\C^\times$ in \cite{[9]}, (Theorem 2, case $n=1$ and its Corollary) and in general in \cite{BGH}.
\subsection{Theorem} \label{ThBGH} The sequence $W^{\frac{1}{m}}$ stops branching if and only if  $W$ is not contained in a coset of a proper algebraic subgroup of  $\Ac^k.$

\bc\lb{Lemma3.2.2.}  Any irreducible closed subset of $V^n$ is of the form $L\cap \ln W,$  
for $W$ Kummer generic in $\exp L.$

 Any closed subset of $V^n$ is analytic in $V^n.$ 

\ec

It is easy now to check that (SI),
 (INT), (CMP), (CC), (AS) and (PS) are satisfied. 

\bc \lb{3.9} The structure $(V,L)$ is  analytic Zariski and 
presmooth. 
\ec

The  reader may notice that the analysis above treats only 
{\em formal} notion of analyticity on the cover $\C$ of  $\C^\times$ but does not address the classical one.
In particular, {\em is the formal analytic decomposition as described by  \ref{Lemma3.2.2.} the same as the actual complex analytic one?} In a private communication
F.Campana answered this question in positive, using a cohomological argument.
M.Gavrilovich proved this and much more general statement in his thesis (see \cite{[8]}, III.1.2) by a similar argument.\\ 

\subsection{Covers in positive characteristic}
Now we look into yet another version of a cover structure which is proven to be
analytic Zariski, a {\bf cover of the one-dimensional algebraic torus over  an algebraically closed field of a positive characteristic.} 

Let  $(V,+)$ be a divisible torsion free abelian group and $\F$ an algebraically closed field of a positive characteristic $p.$ We assume that $V$ and $\F$ are
both of the same uncountable cardinality. Under these assumptions it is easy
to construct a surjective homomorphism  $$\exf: V\to \F^\times.$$ The kernel
of such a homomorphism must be a subgroup which is $p$-divisible but not $q$-divisible for each $q$ coprime with $p.$ One can easily construct $\exf$ so that
$$\ker \exf\cong \Z[\frac{1}{p}],$$
the additive group (which is also a ring) of rationals of the form $\frac{m}{p^n},$ $m,n\in \Z,$ $n\ge 0.$ In fact in this case it is convenient to view $V$ and $\ker \exf$ as 
$\Z[\frac{1}{p}]$-modules.

In this new situation Lemma~\ref{qelim} is still true, with obvious 
alterations, and we can use the definition \ref{Lemma3.2.1} to introduce a topology and the family $L$ as above. The necessary version of Theorem \ref{ThBGH} 
 is proved in \cite{BGH}. Hence
the corresponding versions of \ref{Lemma3.2.2.} follows.

\subsection{Remark} \label{remark} In all the above examples the analytic rank of any nonempty closed subset is 1, that is any closed subset is analytic.


\subsection{$\C_{\exp}$ and other pseudo-analytic structures}

$\C_{\exp},$  the structure  $(\C; +,\cdot, \exp),$ was a prototype of {\bf the field with pseudo-exponentiation} studied by the current author in \cite{[14]}. It was proved (with later corrections, see \cite{BK}) that this structure is quasi-minimal and its (explicitly written) $L_{\omega_1,\omega}(Q)$-axioms are categorical in  all uncountable cardinality. This result has been generalised to many other structures of analytic origin in \cite{BK}, in particular to
the the formal analogue of $\C_{\mathfrak{P}}=(\C; +,\cdot, \mathfrak{P}),$ where $\mathfrak{P}=\mathfrak{P}(\tau, z)$ is the Weierstrass function of variable $z$ with parameter $\tau.$  We call these structures {\bf pseudo-analytic.}

It is a reasonable conjecture to assume that the pseudo-analytic structures of cardinality continuum are isomorphic to their complex prototypes. Nevertheless, even under this conjecture it is 
not known whether   $\C_{\exp},$  $\C_{\mathfrak{P}}$ or any of the other pseudo-analytic structures (which do not satisfy    
\ref{remark}) are analytic Zariski. One may start by defining the family of (formal) closed sets in the structure to coincide with the family of definable subsets which are closed in the metric topology of the complex manifold. The problem then is to conveniently classify such subsets. A suggestion for such a classification may come from the following notion.

   \subsection{Generalised analytic sets}\label{gan1}
In \cite{[5]} we have discussed the following notion of 
generalised analytic subsets of $[\PP^1(\C)]^n$ and, more generally, of 
$[\PP^1(K)]^n$ for $K$ algebraically closed complete valued field.

Let $F\subs \C^2$ be a graph of an entire analytic function and $\bar F$ its
closure in $[\PP^1(\C)]^2.$ It follows from Picard's Theorem that 
$\bar F= F\cup \{ \infty\}\times \PP^1(\C),$ in particular $\bar F$ has
analytic rank
$2.$

{\em Generalised analytic sets} are defined as
the subsets of $[\PP^1(\C)]^n$ for all $n,$ obtained
from classical (algebraic) Zariski closed subsets of $[\PP^1(\C)]^n$ and some number  of  sets of the form
$\bar F$
by applying the positive operations: Cartesian products, finite 
intersections, unions and projections. It is clear by definition that the complex generalised analytic sets are closed 
(but not obvious for the case of $K,$ algebraically closed complete non-Archimedean valued field). 

\subsection{Theorem (see \cite{[5]})}\label{gan2}
 {\em Any generalised analytic set is of finite analytic rank}.

\end{document}